\documentclass[10pt]{amsart}

\newcommand{\ee}{\end{eqnarray}}
\newcommand{\be}{\begin{eqnarray}}

\newcommand{\Diff}{\mbox{\rm Diff}}
\newtheorem{thm}{Theorem}[section]
\newtheorem{lem}[thm]{Lemma}
\newtheorem{prop}[thm]{Proposition}

\newtheorem{conj}{Conjecture}

\title[Closing Lemma]{A $C^{r}$ Closing Lemma for a Class of
Symplectic Diffeomorphisms}

\author{Zhihong Xia \& Hua Zhang}
\address{Department of Mathematics \\Northwestern University \\ Evanston,
  IL 60208}

\thanks{Both authors are supported in part by National Science
Foundation.}

\date{\today}
\email{xia@math.northwestern.edu \& zhang@math.northwestern.edu}
\begin{document}

\begin{abstract}
We prove a $C^r$ closing lemma for a class of partially hyperbolic
symplectic diffeomorphisms. We show that for a generic $C^r$
symplectic diffeomorphism, $r =1, 2, \ldots,$, with two dimensional
center and close to a product map, the set of all periodic points is
dense.
\end{abstract}

\maketitle

\section{Introduction and Main Result}

One of the fundamental problems in dynamical systems is the so-called
$C^r$ closing lemma. The problem goes back to Poincar\'e in his study
of the restricted three body problem. It asks whether periodic points
are dense for a typical symplectic or volume preserving diffeomorphism
on a compact manifold. Let $M$ be a compact manifold, with either a
symplectic or a volume form $\omega$. Let $\Diff_{\omega}^{r}(M)$ be
the set of $C^r$ symplectic or volume-preserving diffeomorphisms on
$M$. A set in a topological space is said to be residual if it is the
intersection of countably many open and dense subsets of of the
topological space. A dynamical property is said to be $C^r$ generic on
$\Diff_{\omega}^{r}(M)$ if there is a residual set $R$ such
that the property holds for every $f \in R$. In the class symplectic
and volume preserving diffeomorphisms, the closing lemma is the
following conjecture:

\begin{conj}[Closing Lemma for symplectic and volume-preserving
    diffeomorphisms] Assume $M$ is compact. There exists a residual
  subset $R \subset \Diff_{\omega}^{r}(M)$ such that if $f \in R$, the
  set of periodic points $P = \{ x \in M \; | \; f^p(x) =x, \mbox{ for
  some integer } p \}$ is dense in $M$.
\end{conj}

Smale \cite{Smale98} listed the problem as one of the mathematical
problems for this century. For $r = 1$, the above conjecture is proved
to be true by Pugh \cite{Pu67} and later improved to various cases by
Pugh \& Robinson \cite{PR83}.  A different proof was given by Liao
\cite{Liao79} and Mai \cite{Mai86}. For higher smoothness $r > 1$,
besides the hyperbolic cases (the Anosov closing lemma, for uniformly
hyperbolic and non-uniformly hyperbolic diffeomorphisms), there is no
non-trivial results. On the other hand, example shows that the local
perturbation method used in the proof of the $C^1$ closing lemma no
longer works for the smoother case. New and global perturbation
methods are required (Gutierrez \cite{Gutierrez87}). M.\ Herman
\cite{He91} has a counter example for the $C^r$ closing lemma with $r$
large for symplectic diffeomorphisms where the symplectic form is not
exact.

In this paper, we prove a $C^r$ closing lemma for arbitrary positive
integer $r$ for a class of partially hyperbolic symplectic
diffeomorphisms.

A diffeomorphism $f: M\rightarrow M$ is partially hyperbolic if the
tangent bundle $TM$ admits a $Tf$ invariant splitting $TM=E^{u}\oplus
E^{c}\oplus E^{s}$ and there is a Riemannian metric on $M$
such that there exist real numbers $ \lambda_{1}> \lambda_{2}> 1>
\mu_{2}> \mu_{1}>0 $
 satisfying
$$m(Tf|_{E^{u}})\geq\lambda_{1}>\lambda_{2}\geq\| Tf|_{E^{c}} \| \geq
m(Tf|_{E_{c}})\geq\mu_{2}> \mu_{1}\geq\|Tf|_{E^{s}}\|>0.$$
Here the co-norm $m(A)$ of a linear operator $A$ between two Banach
spaces is defined by $m(A):= inf_{\|v\|=1} \|A(v)\|=\|A^{-1}\|^{-1}$.

To avoid triviality, we assume at least two of the subbundles are
non-zero. Partial hyperbolicity is a $C^{1}$ open condition, as can be
easily verified by it's associated invariant cone fields.

We remark that our definition of partial hyperbolicity here is not the
most general one. One can allow the parameters
$\lambda_{1},\lambda_{2}, \mu_{1}, \mu_{2}$ to depend on each
trajectory in general cases. The systems that we are considering
satisfy the definition given here.

For symplectic cases, the stable distribution $E^s$ and the unstable
distribution $E^u$ have the same dimension. Moreover, one can choose
the parameters such that $\lambda_1 = \mu_1^{-1}$ and $\lambda_2 =
\mu_2^{-1}$.

We are now ready to state our main theorem.

\begin{thm}
Let $M_1$ be a compact symplectic manifold and $f_1\in \mbox{\rm
Diff}_{\omega_1}^{r}(M_1)$ an Anosov diffeomorphism. Let $M_2$ be a
compact symplectic surface (orientable surface) with an area form
$\omega_2$ and let $f_2 \in \Diff_{\omega_2}^{r}(M_2)$ be an area
preserving diffeomorphism on $M_2$. Let $\omega = \omega_1 + \omega_2$
be the symplectic form defined on $M_1 \times M_2$. Assume that $f_1$
dominates $f_2$, i.e., $f_1 \times f_2: M_1 \times M_2 \rightarrow M_1
\times M_2$ is partially hyperbolic with $TM_2$ as its center
splitting. Then there exists a neighborhood $U$ of $f_1\times f_2$ in
$\Diff_{\omega}^{r}(M_1\times M_2)$ and a residual subset $R \in U$
such that for any $g\in U$, the set of periodic points of $g$ is dense
in $M_1\times M_2$.
\end{thm}

The proof took advantage of the partial hyperbolicity and a recent
result of Xia \cite{Xia04b} on surface diffeomorphisms.

\section{Partial Hyperbolicity and Symplectic Structure}

For a $C^{r}$ partially hyperbolic diffeomorphism $f$ , the stable and
unstable bundles are uniquely integrable and are tangent to foliations
$W^{s}_{f}$ and $W^{u}_{f}$ with $C^{r}$ leaves. $f$ is dynamically
coherent if the distributions $E^{c}$, $E^{c}\oplus E^{s}$ and
$E^{c}\oplus E^{u}$ are integrable, they integrate to foliations
$W^{c}_{f}$, $W^{cs}_{f}$ and $W^{cu}_{f}$ respectively and
$W^{c}_{f}$ and $W^{s}_{f}$ sub-foliate $W^{cs}_{f}$, $W^{c}_{f}$ and
$W^{u}_{f}$ sub-foliate $W^{cu}_{f}$.

\vspace{.3cm}

We have the following proposition from Pugh \& Shub \cite{PS97}.

\begin{prop}
Let $f$ be a partially hyperbolic diffeomorphism. If the center foliation $W^{c}_{f}$ exists and is of class $C^{1}$,
then $f$ is stably dynamically coherent, i.e., any $g$ which is
$C^{1}$ sufficiently close to $f$ is dynamically coherent.
\end{prop}

Let $W$ be a foliation of a compact smooth manifold $M$ whose leaves
are $C^{r}$ immersed submanifolds of dimension $k$. For $x\in M$, we
call a set $P(x)\subset W(x)$ a $C^{r}$ plaque of $W$ at $x$ if $P(x)$
is the image of a $C^{r}$ embedding of the unit ball $B\subset
\mathbb{R}^{k}$ into $W(x)$. A plaquation $\textbf{P}$ for $W$ is a
collection of plaques such that every point $x\in M$ is contained in a
plaque $P\in \textbf{P}$.

Let $f$ be a diffeomorphism such that $W$ is invariant under $f$. A
pseudo orbit $\{x_{n}\}_{n\in \mathbb{Z}}$ respects $\textbf{P}$ if
for each $n$, $f(x_{n})$ and $x_{n+1}$ lie in a common plaque $P\in
\textbf{P}$. $f$ is called plaque expansive with respect to $W$ if
there exists $\epsilon>0$ such that if two $\epsilon$-pseudo orbits
$\{x_{n}\}$ and $\{y_{n}\}$ both respect $\textbf{P}$ and $d(x_{n},
y_{n})< \epsilon$ for all $n$, then $x_{n}$ and $y_{n}$ lie in a
common plaque for all $n$.

Hirsch, Pugh and Shub \cite{HPS77} proved that plaque expansiveness
with respect to the center foliation of a partially hyperbolic
diffeomorphism is a $C^{1}$ open property and is satisfied when we
have smooth center foliation $W^{c}$ (Theorem 7.1 and 7.2 in
\cite{HPS77}).

\vspace{.3cm}

It is clear that under the condition of our main theorem, $f = f_1\times
f_2$ is partially hyperbolic with smooth center foliation, so there
exists neighborhood $U$ of $f_1\times f_2$ in
$\Diff_{\omega}^{r}(M_1\times M_2)$ such that any $g\in U$ is partially
hyperbolic, dynamically coherent and plaque expansive with respect to
the center foliation $W^{c}_{g}$.

\vspace{.3cm}

Ni\c{t}ic\v{a} and T\"{o}r\"{o}k in \cite{NT01} proved the following

\begin{prop}
Let $M$ be a compact manifold with a smooth volume form $\mu$, if
$f\in \Diff^{1}_{\mu}(M)$ is partially hyperbolic, dynamically coherent
and plaque expansive with respect to the center foliation $W^{c}_{f}$,
then the periodic center leaves of $f$ are dense in $M$.
\end{prop}

Now we have the following

\begin{lem}
Under the condition of our main theorem, there exists neighborhood $U$
of $f_1\times f_2$ in $\Diff_{\omega}^{r}(M_1\times M_2)$ such
that any $g\in U$ is partially hyperbolic, dynamically coherent and
the periodic center leaves of $g$ are dense in $M_1\times M_2$.
\end{lem}

The proof is a simple application of the above results and we only
have to note that a symplectic diffeomorphism trivially support an
invariant smooth volume form.

\vspace{.3cm}

The following proposition is also from Ni\c{t}ic\v{a} and
T\"{o}r\"{o}k \cite{NT01}.

\begin{prop}

For a partially hyperbolic diffeomorphism on a compact manifold $M$
with center-stable and center-unstable foliations $W^{cs}$ and
$W^{cu}$, we have the following local product structure property:

There exist constants $\epsilon>0$, $\delta>0$ and $K>1$ such that for
any $x,y\in M$ with $d(x,y)<\epsilon$, the following hold,

1) $W^{s}_{\delta}(x)$ and $W^{cu}_{\delta}(y)$ intersect at a unique
   point $z_{1}$, $W^{u}_{\delta}(x)$ and $W^{cs}_{\delta}(y)$
   intersect at a unique point $z_{2}$, and moreover

$max(d(x,z_{1}), d(y,z_{1}))<Kd(x,y)$, $max(d(x,z_{2}),
d(y,z_{2}))<Kd(x,y)$.

2) $W^{cs}_{\delta}(x)$ and $W^{cu}_{\delta}(y)$ intersect
   transversally, same is true for $W^{cs}_{\delta}(y)$ and
   $W^{cu}_{\delta}(x)$.

3) $W^{cs}_{\delta}(x)\bigcap W^{cu}_{\delta}(x)= W^{c}_{\delta}(x)$
   and $W^{cs}_{\delta}(y)\bigcap W^{cu}_{\delta}(y)=
   W^{c}_{\delta}(y)$.
\end{prop}

Theorem 6.1 of \cite{HPS77} tells us that $\epsilon$, $\delta$ and $K$
are lower semi-continuous with respect to the $C^{1}$ topology on
$\Diff^{1}(M)$.

\vspace{.3cm}

We will need a result for symplectic partially hyperbolic diffeomorphisms.

\begin{lem}
Let $f$ be a symplectic partially hyperbolic diffeomorphism on a
compact symplectic manifold $M$, suppose we have the center foliation
$W^{c}_{f}$, then center manifolds of $f$ are symplectic submanifolds
and the restrictions of $f$ on invariant center leaves are symplectic
diffeomorphisms.
\end{lem}

$\textit{Proof}$. For symplectic partially
hyperbolic diffeomorphism $f$, there exist $\lambda>\mu>1$ such that
$$m(Tf|_{E^{u}})\geq\lambda>\tau\geq\| Tf|_{E^{c}} \| \geq
m(Tf|_{E_{c}})\geq\tau^{-1}> \lambda^{-1}\geq\|Tf|_{E^{s}}\|>0.$$

Denote by $\omega$ the symplectic form on $M$. Let
$W\subset M$ be a center leaf, we should prove $(W,\omega|_{W})$ is a
symplectic manifold, i.e., $\omega|_{W}$ is a non-degenerate, closed
two form on $W$. Closeness is obvious since $\omega$ is closed on $M$.

  Suppose that $\omega|_{W}$ is degenerate, then there exists $x\in W$
  and a unit vector $u\in T_{x}W$ such that for all $v_{c}\in
  T_{x}W$, $\omega(u,v_{c})=0$.

  We have the splitting $T_{x}M=E_{x}^{s}\oplus E_{x}^{c}\oplus
  E_{x}^{u}$, any $v\in T_{x}M$ can be written as
  $v=v_{s}+v_{c}+v_{u}$, where $v_{s}\in E_{x}^{s}$, $v_{u}\in
  E_{x}^{u}$ and $v_{c}\in E_{x}^{c}=T_{x}W$.

  We have  $\omega(u,v)=\omega(u,v_{s})+\omega(u,v_{c})+\omega(u,v_{u})$ and
  by the way $u$ was chosen, $\omega(u,v_{c})=0$.

  There exists $K>0$ such that
  $|\omega(w^{1},w^{2})|<K$ for arbitrary pair of unit vectors
  $w^{1},w^{2}\in T_{x}M$ .

  Now we know $\omega(u,v_{s})=0$ because 

  $|\omega(u,v_{s})|=|\omega(Tf^{n}(u),Tf^{n}(v_{s})|
  \leq(\frac{\tau}{\lambda})^{n}\|v_{s}\||\omega(u^{n},v^{n}_{s})|\leq
  K(\frac{\tau}{\lambda})^{n}\|v_{s}\|\longrightarrow0$ as
  $n\longrightarrow +\infty$.

  Here $u^{n}=Tf^{n}(u)/\|Tf^{n}(u)\|$,
  $v^{n}_{s}=Tf^{n}(v_{s})/\|Tf^{n}(v_{s})\|$.

  Similarly, we have $\omega(u,v_{u})=0$ and hence $\omega(u,v)=0$ for
  any $v\in T_{x}M$, this contradicts the fact that $\omega$ is
  non-degenerate on $M$.

  So $(W,\omega|_{W})$ is a symplectic submanifold and if $W$ is
  invariant under $f$, $f|_{W}$ preserves $\omega|_{W}$ and hence is a
  symplectic diffeomorphism on $W$.

\section{Some generic properties for area-preserving diffeomorphisms
  on surfaces}

To prove our main theorem, we need some generic properties for surface
diffeomorphisms. Let $S$ be a compact surface, denote by
$\Diff^{r}_{\mu}(S)$ the set of area-preserving $C^{r}$
diffeomorphisms. For $f\in \Diff^{r}_{\mu}(S)$, denote by $HP(f)$ the
set of hyperbolic periodic points of $f$. The following generic
property was first proved by Mather \cite{Mather82a} for maps on two
sphere $S^2$ and later generalized to arbitrary compact surfaces by
Oliveira \cite{Oliveira87}.

\begin{prop}
There is a residual subset $R \in \Diff^r_\mu(S)$ such that if $f \in
R$ and $p \in HP(f)$ is a hyperbolic periodic point of $f$, then
$$\overline{W^{s}_{f}(p)} = \overline{W^{u}_{f}(p)}.$$
\end{prop}

We remark that if $r=1$, the proposition is true for generic
symplectic and volume preserving diffeomorphisms on any compact
manifolds (cf. Xia \cite{Xi96a}).

The next Theorem is due to Xia \cite{Xia04b}, extending a recent
result of Franks \& Le Calvez \cite{FL03} on two sphere.

\begin{thm}
Let $S$ be a compact orientable surface and $\mu$ be an area form on
$S$. For any positive integer $r$, there is a residual subset $R \in
\Diff^r_\mu(S)$ such that if $f \in R$, then both the sets $\cup_{p\in
HP(f)} W^s(p)$ and $\cup_{p\in HP(f)} W^u(p)$ are dense in
$S$. Moreover, if an open set $U \subset S$ contains no periodic
point, then there is a hyperbolic periodic point $p \in HP(f)$ such
that both the stable and unstable manifold of $p$ is dense in $U$.
\end{thm}

The proof uses prime end compactification and a rich literature on
area preserving surface diffeomorphisms.

\section{Proof of the Main Theorem}

Our main perturbation lemma is from Xia \cite{Xi96a}.

\begin{lem}
Let $M$ be a compact symplectic manifold and $f\in
\Diff^{r}_{\omega}(M)$, there exist $\epsilon_{0}>0$ and $c>0$
such that for any $g\in \Diff^{r}_{\omega}(M)$ such that
$\|f-g\|_{C^{r}}<\epsilon_{0}$ and any $0<\epsilon\leq
\epsilon_{0}$, $0<\delta\leq \epsilon_{0}$, if $x, y\in M$ and
$d(x,y)<c\delta^{r}\epsilon$, there exists $g_{1}\in
\Diff^{r}_{\omega}(M)$, $\|g_{1}-g\|_{C^{r}}<\epsilon$ satisfies
$g_{1}g^{-1}(x)=y$, $g_{1}(z)=g(z)$ for all $z\notin
g^{-1}(B_{\delta}(x))$, and $g_{1}^{-1}(z)=g^{-1}(z)$ for all $z\notin
B_{\delta}(x)$.

\end{lem}

Now we are ready to prove the main theorem.

\vspace{.3cm}

$\textit{Proof}$. By Lemma 2.3, there exists a neighborhood $U$ of
$f_1\times f_2$ in $\Diff_{\omega}^{r}(M_1\times M_2)$ such that
any $g\in U$ is partially hyperbolic, dynamically coherent and the
periodic center leaves of $g$ are dense in $M_1\times M_2$.

Now for any fixed $g\in U$, suppose there is a periodic point free
open set $V\subset M_1\times M_2$, we show that by an arbitrarily small
perturbation, we can create a periodic point in $V$. It is clear that
the main theorem will follow.

Since periodic center leaves of $g$ are dense, by Proposition 2.4, we
can find two periodic center leaves $W_{1}$ and $W_{2}$ which are
sufficiently close such that there exist $x_{1},y_{1}\in W_{1}\bigcap
V$, $x_{2},y_{2}\in W_{2}\bigcap V$ with
$z=W^{u}_{\delta}(x_{1})\bigcap W^{s}_{\delta}(x_{2})\in V$,
$w=W^{s}_{\delta}(y_{1})\bigcap W^{u}_{\delta}(y_{2})\in V$.  $W_{1}$ and $W_{2}$ are compact surfaces.

By taking certain power of $g$ we may assume that $W_{1}$ and $W_{2}$
are invariant under $g$. From Lemma 2.5, $W_{1}$ and $W_{2}$ are
symplectic submanifolds, $g^{1}=g|_{W_{1}}$ and $g^{2}=g|_{W_{2}}$ are
symplectic diffeomorphisms. By making an arbitrarily small
perturbation, we may assume that $g^{1}$ and $g^{2}$ satisfy the
generic condition in Theorem 3.2. Now $W_{i}\bigcap V$ is a periodic
point free open set in $W_{i}$, we know that there exists $p_{i}\in
HP(g^{i})$ such that $W^{u}_{g^{i}}(p_{i})$ and $W^{s}_{g^{i}}(p_{i})$
are both dense in $W_{i}\bigcap V$, where $i=1,2$. Note that $p_{1}$
and $p_{2}$ are hyperbolic periodic points of $g$.

We will show that by an arbitrarily small perturbation, we can change
$z$ and $w$ into heteroclinic points of hyperbolic periodic points
$p_{1}$ and $p_{2}$ and get a heteroclinic loop. As a result, there
will be periodic points in arbitrary neighborhoods of $z$ and $w$,
including $V$.

For arbitrary $\eta>0$ prescribed as the size of the perturbation,
take $\epsilon$ such that $0<\epsilon<
min(\frac{\eta}{4},\epsilon_{0})$, where $\epsilon_{0}$ is from Lemma
4.1.

Since $z\in W^{u}(W_{1})$ and $z\notin W_{1}$, there exists $\alpha$
with $0<\alpha<\epsilon_{0}$, such that $B_{\alpha}(z)\bigcap
W_{1}=\emptyset$ and $g^{-n}(z)\notin B_{\alpha}(z)$ for all $n\in
\mathbb{N}$.

Fix this $\alpha$, there exists $\beta$ with $0<\beta<\alpha$ such
that for all $\tilde{z}\in W^{u}(W_{1})$ with $d(\tilde{z}, z)<\beta$,
we have $g^{-n}(\tilde{z})\notin B_{\frac{\alpha}{2}}(z)$ for all
$n\in \mathbb{N}$.

Fix this $\beta$, there exists $\gamma$ with $0< \gamma< min(\beta,
c(\frac{\alpha}{4})^{r}\epsilon)$, such that for all $z_{1}\in
W^{u}(W_{1})$ with $d(z_{1}, z)<\gamma$, we have
$B_{\frac{\alpha}{4}}(z_{1})\subset B_{\frac{\alpha}{2}}(z)$ and hence
$g^{-n}(z_{1})\notin B_{\frac{\alpha}{4}}(z_{1})$ for all $n\in
\mathbb{N}$.

By continuity of the unstable foliation, there exists $\nu>0$ such
that for all $\tilde{x}_{1}\in W_{1}$ with $d(\tilde{x}_{1},
x_{1})<\nu$, there exists $z_{1}\in W^{u}_{g}(\tilde{x}_{1})$ such
that $d(z_{1}, z)<\gamma$.

Since $W^{u}_{g^{1}}(p_{1})$is dense in $W_{1}\bigcap V$, there exists
$\tilde{x}_{1}\in W^{u}_{g^{1}}(p_{1})$ with $d(\tilde{x}_{1},
x_{1})<\nu$ and hence there is a $z_{1}\in W^{u}_{g}(\tilde{x}_{1})$
such that $d(z_{1}, z)<\gamma$.

Now we can use the perturbation lemma 4.1 for $g$ using the parameters
$0<\epsilon<\epsilon_{0}$ and
$0<\delta=\frac{\alpha}{4}<\epsilon_{0}$. We have $d(z_{1},
z)<\gamma<c\delta^{r}\epsilon$, so there exists $g_{1}\in
\Diff^{r}_{\omega}(M_1\times M_2)$, $\|g_{1}-g\|_{C^{r}}<\epsilon$
satisfies $g_{1}g^{-1}(z_{1})=z$, $g_{1}(x)=g(x)$ for all $x\notin
g^{-1}(B_{\delta}(z_{1}))$, and $g_{1}^{-1}(y)=g^{-1}(y)$ for all
$y\notin B_{\delta}(z_{1})$.

We check that after the perturbation, $z\in \tilde{W}^{u}_{g_{1}}(p_{1})$, where 
$\tilde{W}^{u}_{g_{1}}(p_{1})$ stands for the unstable manifold of the hyperbolic periodic point $p_1$ for $g_1$, not the leaf of the unstable foliation containing $p_1$ in the partially hyperbolic setting.
 
 It is clear that $g_{1}^{-1}(z)=g^{-1}(z_{1})$ and since
$g^{-n}(z_{1})\notin B_{\delta}(z_{1})$ for all $n\in \mathbb{N}$,
$g_{1}^{-n}(z)=g^{-n}(z_{1})$ for all $n\in \mathbb{N}$. Moreover,
$g_{1}^{-n}(p_{1})=g^{-n}(p_{1})$ for all $n\in \mathbb{N}$ since
$B_{\delta}(z_{1})\bigcap W_{1}=\emptyset$.

Hence we have 

$d(g_{1}^{-n}(z),g_{1}^{-n}(p_{1}))=d(g^{-n}(z_{1}),g^{-n}(p_{1}))$

$\leq d(g^{-n}(z_{1}), g^{-n}(\tilde{x}_{1})) +
d(g^{-n}(\tilde{x}_{1}), g^{-n}(p_{1})) \longrightarrow 0 $ as
$n\longrightarrow +\infty$.

This shows $z\in \tilde{W}^{u}_{g_{1}}(p_{1})$. The two terms above both go to
$0$ as $n$ goes to $+\infty$ since $z_{1}\in W^{u}_{g}(\tilde{x}_{1})$
and $\tilde{x}_{1}\in W^{u}_{g^{1}}(p_{1})$.

Similarly we can use a perturbation of size less than $\epsilon$ to
make $z$ on the stable manifold of $p_{2}\in W_{2}$. Two more of these
perturbations will make $w$ in the intersection of stable manifold of
$p_{1}$ and unstable manifold of $p_{2}$. Finally by a perturbation of
size less than $4\epsilon=\eta$, we have the desired heteroclinic
loop. This concludes our proof.


\end{document}